\documentclass[english]{elsarticle}
\usepackage[T1]{fontenc}
\usepackage[utf8]{inputenc}
\usepackage[a4paper]{geometry}
\usepackage{babel}
\usepackage{url}
\usepackage{amsmath}
\usepackage{amssymb}
\usepackage{graphicx}
\usepackage[unicode=true,
 bookmarks=true,bookmarksnumbered=false,bookmarksopen=false,
 breaklinks=false,pdfborder={0 0 1},backref=false,colorlinks=false]
 {hyperref}

\makeatletter
\usepackage{tikz}

\makeatother

\begin{document}

\begin{frontmatter}{}

\title{Palindromic discontinuous Galerkin method for kinetic equations with
stiff relaxation}

\author{David Coulette, Emmanuel Franck, Philippe Helluy, Michel Mehrenberger,
Laurent Navoret}

\ead{helluy@unistra.fr}

\address{IRMA, Univ. Strasbourg, 7 rue Descartes, Strasbourg, France \& Inria
TONUS}
\begin{abstract}
We present a high order scheme for approximating kinetic equations
with stiff relaxation. The objective is to provide efficient methods
for solving the underlying system of conservation laws. The construction
is based on several ingredients: (i) a high order implicit upwind
Discontinuous Galerkin approximation of the kinetic equations with
easy-to-solve triangular linear systems; (ii) a second order asymptotic-preserving
time integration based on symmetry arguments; (iii) a palindromic
composition of the second order method for achieving higher orders
in time. The method is then tested at orders 2, 4 and 6. It is asymptotic-preserving
with respect to the stiff relaxation and accepts high CFL numbers.
\end{abstract}

\end{frontmatter}{}

\section{Introduction}

The Lattice Boltzmann Method (LBM) is a general method for solving
systems of conservation laws \cite{chen1998lattice}. Today, it is
routinely applied in fluid flow simulations. The LBM relies on a kinetic
representation of the system of conservation laws by a small set of
transport equations coupled through a stiff relaxation source term.
The transport velocities are generally taken at the vertices of a
regular lattice and the construction resembles the Boltzmann kinetic
theory of gases, hence the name of the method. The kinetic model is
solved with a splitting method, in which the transport and relaxation
steps are treated separately. Usually, the transport is exactly solved
by the characteristic method, while the relaxation is approximated
by a second order Crank-Nicolson scheme \cite{dellar2013interpretation}. 

The main drawback of the LBM is that it requires regular grids and
that the time step $\Delta t$ is imposed by the grid step $\Delta x$.
In addition, with some of the most widely used lattices, the stability
properties of the method are sometimes unclear \cite{dubois2013stable,helluy:hal-01403759}.

Several authors have proposed methods for relaxing the grid constraint.
It is possible to solve the kinetic model with a Finite Difference
scheme (FDLBM) \cite{mei1998finite}, a Finite Volume method (FVLBM)
\cite{nannelli1992lattice,peng1998lattice} or a Discontinuous Galerkin
approximation (DGLBM) \cite{jin1999efficient}.

Recently, Graille \cite{graille2014approximation} has proposed a
vectorial LBM that he was able to reinterpret as a relaxation system
in the sense of Jin and Xin \cite{jin1995relaxation} with proven
stability properties.

In this paper, we propose a DGLBM based on the Graille's vectorial
LBM.

We extend the original DGLBM (see \cite{shi2003discontinuous,min2011spectral})
in several directions. The first improvement is to apply an implicit
DG method instead of an explicit one for solving the transport equations.
This can be done at almost no additional cost. Indeed, with an upwind
numerical flux, the linear system of the implicit DG method is triangular
and, in the end, can be solved explicitly. In this way, we obtain
stable methods even with high CFL numbers. This kind of ideas can
be found for instance in \cite{bey1997downwind,wang1999crosswind,coquel2008large,natvig2008fast,moustafa20153d}.

The second improvement is to construct a symmetric-in-time integrator
that remains second order accurate even for vanishing relaxation time
(asymptotic-preserving property \cite{jin1999efficient}).

The third improvement is then to embed the second order scheme into
a higher order method thanks to a palindromic composition (described
for instance in \cite{suzuki1990fractal,kahan1997composition,mclachlan2002splitting,hairer2006geometric}).
We have tested fourth and sixth order time integration.

The objective of this paper is to present the whole construction of
the Palindromic Discontinuous Galerkin Lattice Boltzmann Method and
then to apply it on one-dimensional test cases.

\section{A general lattice kinetic interpretation of systems of conservation
laws }

\subsection{One-dimensional case}

We are interested in the numerical simulation of a system of conservation
laws
\begin{equation}
\partial_{t}w+\partial_{x}q(w)=0.\label{eq:conslaw}
\end{equation}
The unknown is a vector valued function $w(x,t)\in\mathbb{R}^{m}$
depending on a space variable $x\in[a,b]$ and time $t\in[0,T]$.
The flux $q(w)\in\mathbb{R}^{m}$ is such that (\ref{eq:conslaw})
is a hyperbolic system of conservation laws: its Jacobian matrix $D_{w}q(w)$
is diagonalizable with real eigenvalues.

Using the approach of \cite{graille2014approximation} it is possible
to construct a general kinetic interpretation of (\ref{eq:conslaw})
in the following way: we consider a family of $2m$ transport equations
with a relaxation source term
\begin{equation}
\partial_{t}f_{i}+v_{i}\partial_{x}f_{i}=\frac{1}{\tau}(f_{i}^{eq}-f_{i}),\label{eq:kinetic}
\end{equation}
where the microscopic distribution function $f_{i}(x,t)$ depends
on the space variable $x$, a velocity index $i\in\{0,\ldots,2m-1\}$,
and time $t$. 

The velocities $v_{i}$ are taken in a lattice with $2m$ points:
\begin{equation}
v=(v_{0},v_{1},\ldots,v_{2m-2},v_{2m-1})=(-\lambda,\lambda,\ldots,-\lambda,\lambda),\label{eq:velocity_lattice}
\end{equation}
where $\lambda$ is a positive constant.

The conservative variables $w=(w_{0}\ldots w_{m-1})$ are related
to $f$ through
\[
w_{k}=\sum_{i=0,1}f_{2k+i}.
\]
We also define $z$, the vector of first moments of $f$
\[
z_{k}=\sum_{i=0,1}v_{i}f_{2k+i}.
\]
If we introduce the matrices
\[
M_{w}=\left(\begin{array}{ccccc}
1 & 1 & \cdots & 0 & 0\\
\vdots & \vdots & \ddots & \vdots & \vdots\\
0 & 0 & \cdots & 1 & 1
\end{array}\right),\quad M_{q}=\left(\begin{array}{ccccc}
-\lambda & \lambda & \cdots & 0 & 0\\
\vdots & \vdots & \ddots & \vdots & \vdots\\
0 & 0 & \cdots & -\lambda & \lambda
\end{array}\right),\quad M=\left(\begin{array}{c}
M_{w}\\
M_{q}
\end{array}\right),
\]
the above relations can also be written in matrix form
\begin{equation}
Mf=\left(\begin{array}{c}
w\\
z
\end{array}\right).\label{eq:out_of_eq}
\end{equation}

The equilibrium distribution function $f^{eq}(f)$ is defined in such
a way that 
\begin{eqnarray}
w_{k} & = & \sum_{i=0,1}f_{2k+i}^{eq},\quad q(w)_{k}=\sum_{i=0,1}v_{2k+i}f_{2k+i}^{eq}.\label{eq:w_cons}
\end{eqnarray}
This gives
\begin{equation}
f_{2k}^{eq}=\frac{w_{k}}{2}-\frac{q(w)_{k}}{2\lambda},\quad f_{2k+1}^{eq}=\frac{w_{k}}{2}+\frac{q(w)_{k}}{2\lambda}.\label{eq:maxwellian}
\end{equation}
The above relations can also be written in matrix form
\begin{equation}
Mf^{eq}=\left(\begin{array}{c}
w\\
q(w)
\end{array}\right).\label{eq:in_eq}
\end{equation}

Introducing the diagonal matrix
\[
V=\left(\begin{array}{cccc}
v_{0} &  &  & 0\\
 & v_{1}\\
 &  & \ddots\\
0 &  &  & v_{2m-1}
\end{array}\right)
\]
the kinetic equations (\ref{eq:kinetic}) become
\begin{equation}
\partial_{t}f+V\partial_{x}f=\frac{1}{\tau}(f^{eq}(f)-f).\label{eq:vect_kin}
\end{equation}
Multiplying (\ref{eq:vect_kin}) on the left by $M$ and using (\ref{eq:out_of_eq}),
(\ref{eq:in_eq}) we get
\begin{align*}
\partial_{t}w+\partial_{x}z & =0,\\
\partial_{t}z+\lambda^{2}\partial_{x}w & =\frac{1}{\tau}(q(w)-z).
\end{align*}

As explained in \cite{graille2014approximation}, this shows that
the kinetic model is equivalent to a relaxation approximation of (\ref{eq:conslaw}).
When $\tau\to0$, then $z=q(w)$ and $w$ formally satisfies the initial
system (\ref{eq:conslaw}).

The model (\ref{eq:kinetic}), (\ref{eq:out_of_eq}), (\ref{eq:maxwellian})
is thus a minimalist abstract kinetic interpretation of any system
of conservation laws. By analogy with the Boltzmann theory of gases,
the equilibrium function $f^{eq}(f)$ may be called a ``Maxwellian''
and the relaxation source term the ``BGK'' or ``collision'' term.

The kinetic approximation is stable under the sub-characteristic condition
\[
\lambda>\text{ spectral radius of }D_{w}q(w)),
\]
which states that the lattice velocities must be faster than the fastest
wave of the approximated hyperbolic system \cite{witham1974linear,liu1987hyperbolic,jin1995relaxation}.

\subsection{Generalization to higher dimensions}

In this paper, we shall only consider numerical applications in the
one-dimensional case $x\in\mathbb{R}^{D}$ with $D=1$, but the approach
could be extended to higher dimensions $D>1$. For $D=1$ we have
considered $D+1=2$ kinetic velocities $(\lambda_{0},\lambda_{1})=(-\lambda,\lambda)$.
In higher dimensions $x=(x_{1}\ldots x_{D})$, the system of conservation
laws and the kinetic equations read respectively
\[
\partial_{t}w+\sum_{j=1}^{D}\partial_{x_{j}}q^{j}(w)=0,\quad\partial_{t}f_{i}+\sum_{j=1}^{D}v_{i}^{j}\partial_{x_{j}}f_{i}=\frac{1}{\tau}(f_{i}^{eq}-f_{i}).
\]

The kinetic velocity $\lambda_{i}=(\lambda_{i}^{1}\ldots\lambda_{i}^{D})\in\mathbb{R}^{D}$,
$i=0\ldots D$, can be taken at the vertices of a regular $D-$simplex
for instance. The lattice velocities are then defined by
\[
v_{(D+1)k+i}=\lambda_{i},\quad k=0\ldots m-1,\quad i=0\ldots D.
\]
The equations satisfied by the equilibrium distribution function become
\begin{eqnarray}
w_{k} & = & \sum_{i=0}^{D}f_{(D+1)k+i}^{eq},\quad q^{j}(w)_{k}=\sum_{i=0}^{D}v_{i}^{j}f_{(D+1)k+i}^{eq},\quad j=1\ldots D.\label{eq:w_cons_multi_d}
\end{eqnarray}

Programming and testing the method in two or three dimensions is the
objective of a forthcoming work.

\section{First order CFL-less approximation}

Let us now consider an approximation $f_{h}(t)$ of $f(\cdot,t)$
in a finite-dimensional space $E_{h}$. We assume that the approximation
error behaves like $O(h^{p})$ with $p\geq1$: the space approximation
is at least first order accurate with respect to the discretization
parameter $h$. The kinetic equation (\ref{eq:kinetic}) is thus approximated
by a set of differential equations
\begin{equation}
\partial_{t}f_{h}+L_{h}f_{h}+N_{h}f_{h}=0,\label{eq:kinetic_approx}
\end{equation}
where the operator $L_{h}$ is an approximation of the transport operator
$L=v\partial_{x}$ and $N_{h}$ an approximation of the BGK operator
$N$. Many possibilities may be envisaged: finite differences, finite
elements, discrete Fourier transform, Discontinous Galerkin (DG) approximation,
semi-Lagrangian methods, \textit{etc}. Even if $v\partial_{x}$ is
a linear operator $L_{h}$ might be non-linear (if a limiter technique
is applied in a finite volume method for instance). In this paper,
we adopt an upwind DG approximation with Gauss-Lobatto interpolation
\cite{hesthaven2007nodal}. The collision BGK operator $Nf=-1/\tau(f^{\text{eq}}-f)$
is approximated by the discrete collision operator
\[
N_{h}f_{h}=-\frac{1}{\tau}(f_{h}^{\text{eq}}-f_{h}),
\]
where $f_{h}^{\text{eq}}$ is a discrete equilibrium state, computed
from $f_{h}$ by applying formula (\ref{eq:maxwellian}). It is non-linear
because $f\mapsto f^{eq}$ is non-linear.

The exact flow of the differential equation (\ref{eq:kinetic_approx})
is given by
\[
f_{h}(t)=\exp(-t(L_{h}+N_{h}))f_{h}(0).
\]

The exponential notation can be made completely rigorous here even
in the case of non-linear operators thanks to the Lie algebra formalism.
For an exposition of this formalism in the context of numerical methods
for ordinary differential equations, we refer for instance to the
book \cite{hairer2006geometric} or to the review \cite{mclachlan2002splitting}.

Computing the exact flow is generally not possible. We propose to
apply a splitting method in order to integrate the differential equation
(\ref{eq:kinetic_approx}).

We consider approximations of the collision and transport exact time
integrators that are first order in time. For instance, $C_{1}$ and
$T_{1}$ can be obtained from the implicit first order Euler scheme
\[
C_{1}(\Delta t)=(Id+\Delta tN_{h})^{-1},\quad T_{1}(\Delta t)=(Id+\Delta tL_{h})^{-1}.
\]
Let us point out that $C_{1}$ is a non-linear operator, because $f\mapsto f^{eq}$
is non-linear. The linearity of $T_{1}$ depends on the linearity
of $L_{h}$.

For a fixed $\tau>0$, we have the estimates\footnote{For one single time step the error is $O(\Delta t^{2})$. But when
the error is accumulated on $t_{\max}/\Delta t$ time steps it indeed
produces a first order method.}
\[
C_{1}(\Delta t)=\exp(-\Delta tN_{h})+O(\Delta t^{2}),\quad T_{1}(\Delta t)=\exp(-\Delta tL_{h})+O(\Delta t^{2}).
\]

Finally, let us point out that even if $C_{1}$ and $T_{1}$ are implicit
operators, they are actually very easy to compute. Indeed, $w$ and
thus $f^{eq}$ are constant during the collision step and then $C_{1}$
is simply given by
\[
C_{1}(\Delta t)f=\frac{f^{eq}(f)+\frac{\tau}{\Delta t}f}{1+\frac{\tau}{\Delta t}}.
\]
In addition, because the free transport step is solved by an upwind
DG solver, then the linear operator $Id+\Delta tL_{h}$ is triangular
and its inverse can also be computed explicitly. In practice the transport
step can be solved by scanning the mesh from left to right when $v_{i}>0$
or from right to left when $v_{i}<0$. It is also possible to assemble
the sparse matrix of the implicit DG solver. Some numerical linear
solvers are then able to detect that the matrix is triangular and
apply adequate optimization. This is the case of the KLU library \cite{davis2010algorithm}.
Further optimization can be achieved if we observe that the implicit
solver does not require to keep the values of $f$ at the previous
time step. They can be overwritten by the new values as they are computed.
In general this storage saving is not possible with an explicit DG
scheme applied to a general hyperbolic system.

We can consider the simple Lie's splitting method 
\[
M_{1}=C_{1}T_{1}
\]
for constructing a first order time integrator of the differential
equation (\ref{eq:kinetic_approx}). For a single time-step we have
the following estimate 
\begin{equation}
f_{h}(\Delta t)=\exp\left(-\Delta t(L_{h}+N_{h})\right)f_{h}(0)=C_{1}(\Delta t)T_{1}(\Delta t)f_{h}(0)+O(\Delta t^{2}).\label{eq:liesplitting}
\end{equation}

The resulting first order method $M_{1}$ enjoys good properties: 
\begin{itemize}
\item it is unconditionally stable in time because of the underlying implicit
steps;
\item it is uniformly first order in $\Delta t$, independently of $\tau$
(``asymptotic-preserving'' property \cite{jin1999efficient}) if
the initial condition is sufficiently close to the equilibrium manifold
$\{f,\;f=f^{eq}(f)\}$;
\item it permits low storage optimization;
\item despite being formally implicit, in practice it only requires explicit
computations;
\item it is well-adapted to parallel optimization (and in higher dimensions
$D>1$ the optimization opportunities are even better \cite{moustafa20153d}).
\end{itemize}
However it suffers from an insufficient precision in time.

In this paper we apply techniques that are well-known by the community
of geometric integration for improving the order of (\ref{eq:liesplitting}).
The difficulty is, of course, to keep high order even when $\tau\to0$.
We consider a numerical scheme for computing $\exp(-t(L_{h}+N_{h}))$
that is based on a improved symmetric Strang-splitting scheme associated
to composition methods described in \cite{suzuki1990fractal,kahan1997composition,mclachlan2002splitting,castella2009splitting}
for achieving high order. We show numerical results confirming the
excellent precision and stability of the method.

\section{Second order symmetric time-stepping}

\subsection{Symmetric method}

The first ingredient is to improve the time order of the transport
and collision steps. This can be achieved by symmetrizing the integrator.
For a first order integrator $S_{1}=C_{1}$ or $S_{1}=T_{1}$ we can
construct another integrator by the symmetrization formula
\[
S_{2}(\Delta t)=S_{1}(-\frac{\Delta t}{2})^{-1}S_{1}(\frac{\Delta t}{2}).
\]
The operator $S_{2}$ satisfies the time symmetry property
\begin{equation}
S_{2}(-\Delta t)=S_{2}(\Delta t)^{-1},\label{eq:sym_prop}
\end{equation}
which expresses that if we apply the same method with the opposite
time to the final solution, then we recover exactly the initial condition.
This construction is well known (see for instance \cite{hairer2006geometric,mclachlan2002splitting}).
Actually, formula (\ref{eq:sym_prop}) is an alternative abstract
way to construct the trapezoidal method. Because of symmetry, it is
necessary of order 2 in $\Delta t$ for a fixed $\tau>0$.

The second order collision step is then given by
\begin{equation}
C_{2}(\Delta t)f=\frac{(2\tau-\Delta t)f}{2\tau+\Delta t}+\frac{2\Delta tf^{eq}}{2\tau+\Delta t}.\label{eq:crank_nicolson}
\end{equation}

The second order transport operator $T_{2}$ is constructed in the
same way
\[
T_{2}(\Delta t)=(Id-\frac{\Delta t}{2}L_{h})(Id+\frac{\Delta t}{2}L_{h})^{-1}.
\]
It preserves the nice features of the first order operator $T_{1}$:
stability, parallelism, low storage, triangular matrix structure.

From the two second order bricks $T_{2}$ and $C_{2}$ we can now
construct a global second-order method using Strang formula
\begin{equation}
M_{2}(\Delta t)=T_{2}(\frac{\Delta t}{2})C_{2}(\Delta t)T_{2}(\frac{\Delta t}{2}).\label{eq:bad_m2}
\end{equation}
Because $M_{2}$ satisfies the symmetry property (\ref{eq:sym_prop})
it is necessary of order 2 for a fixed relaxation time $\tau>0$.

\textbf{Remark}: a popular method for solving kinetic equations is
the Lattice Boltzmann Method (LBM). In the LBM the transport step
is solved with an exact characteristic formula. This is possible on
a Cartesian grid when the time step $\Delta t$, the grid step $\Delta x$
and the velocities $v_{i}$ are chosen properly. The collision differential
equation can be solved by several different methods. However, its
stiffness when $\tau$ is small requires implicit or exact exponential
methods. The Euler implicit scheme gives only a first order approximation,
which is not surprising. More surprisingly, the exact exponential
integration also gives poorly accurate results. It is well known in
the LBM community that it is better to consider a trapezoidal second
order time integration for the collision step \cite{dellar2013interpretation}.
This is exactly what we do in our method too. 

\subsection{Asymptotic-preserving second order method }

The above construction works perfectly when the relaxation time $\tau$
is not too small. However, one observes a loss of precision when $\tau\to0$
and the method then returns to first order. This order reduction is
well-known and several cures have been proposed (see for instance
\cite{jin1995runge} and related works).

The phenomenon can be explained in a very simple way. Indeed, if the
method were an approximation of the exact integrator,
\[
M_{2}(t)\simeq\exp(-t(L_{h}+N_{h})),
\]
then we should observe that
\[
M_{2}(0)=Id.
\]
But here, when $\tau\to0$, the trapezoidal collision integrator becomes
\begin{equation}
C_{2}(\Delta t)f=2f^{eq}-f,\label{eq:weak-project}
\end{equation}
and then
\begin{equation}
M_{2}(0)f=2f^{eq}-f.\label{eq:tau0_sym}
\end{equation}
If we restrict $M_{2}$ to the equilibrium manifold $\left\{ f,\;f=f^{eq}(f)\right\} $
we have $M_{2}(0)f=f$ but it is no more true for arbitrary $f$ close
to the equilibrium manifold. A simple way to recover second order
is to apply anyway the method and to perform a final projection of
the result on the equilibrium manifold \cite{dellar2013interpretation}.
This trick, however, destroys the symmetry property.

In this paper we rather propose the following second order method
\begin{equation}
M_{2}(\Delta t)=T_{2}(\frac{\Delta t}{4})C_{2}(\frac{\Delta t}{2})T_{2}(\frac{\Delta t}{2})C_{2}(\frac{\Delta t}{2})T_{2}(\frac{\Delta t}{4}),\label{eq:good_m2}
\end{equation}
which preserves the symmetry property and allows to recover
\[
M_{2}(0)f=f
\]
even when $\tau=0$.

\section{Higher order extension with compositions\label{sec:Higher-order-extension}}

\subsection{The composition method}

From the elementary symmetric brick it is classic to construct higher
order methods by the composition method, which we now recall.

Let $M_{p}$ be a symmetric method of order $p\geq2$. By symmetry,
$p$ is necessary even. A simple approach is to look for a higher
order symmetric method under the form
\[
M_{p+2}(t)=M_{p}(\alpha_{p}t)M_{p}(\beta_{p}t)M_{p}(\alpha_{p}t).
\]
The new method is of order $p+2$ if 
\begin{equation}
2\alpha_{p}+\beta_{p}=1,\quad2\alpha_{p}^{p+1}+\beta_{p}^{p+1}=0.\label{eq:order_cond}
\end{equation}
The proof is elementary (see for instance \cite{mclachlan2002splitting}).
This set of equations admits one pair of real solutions 
\begin{equation}
\alpha_{p}=\frac{1}{2-2^{1/(p+1)}},\quad\beta_{p}=-\frac{2^{1/(p+1)}}{2-2^{1/(p+1)}}.\label{eq:real_coefficients}
\end{equation}
Unfortunately, $\beta_{p}<0$ and the resulting method leads to applying
$M_{p}$ with negative time steps. This may give poor stability properties
when directly applied to dissipative operators. For instance, the
use of a negative time step $\Delta t$ in $C_{2}$ can be catastrophic.
From (\ref{eq:crank_nicolson}) we see that we get a division by zero
when $\Delta t=-2\tau$. In principle, it should not be a problem
when $\tau=0$ because from (\ref{eq:tau0_sym}) we see that the discrete
collision operator is perfectly reversible (actually, it does not
depend on $\Delta t$ anymore).

\subsection{Complex time steps}

In our situation we can also consider complex solutions of (\ref{eq:order_cond}).
A possible choice is to take
\begin{equation}
\alpha_{p}=\frac{\exp(\frac{i\pi}{p+1})}{2\exp(\frac{i\pi}{p+1})+2^{\frac{1}{p+1}}},\quad\beta_{p}=\frac{2^{\frac{1}{p+1}}}{2\exp(\frac{i\pi}{p+1})+2^{\frac{1}{p+1}}}.\label{eq:complex_coefficients}
\end{equation}
This choice requires to apply the basic methods $C_{2}(\Delta t)$
and $T_{2}(\Delta t)$ with a complex time step $\Delta t$. Compared
to the real composition method (\ref{eq:real_coefficients}) the real
part of $\Delta t$ is now $>0$. This can improve in some cases the
stability of the whole composition.

In practice the complex time steps pose no particular problem. In
our C99 implementation, it was enough to replace the ``double''
declarations by ``double complex'' numbers at the adequate places.
In some applications, however, the extension to complex numbers could
be more difficult. This would be the case for instance if the model
requires to using non-analytic boundary conditions or if the transport
solver relies on non-analytic slope limiters. 

For a more detailed analysis and description of the complex composition
method we refer to \cite{castella2009splitting} and included references.
In \cite{castella2009splitting} it is shown that the approach can
be extended up to order $p=14$ (with an additional trick). For $p\geq16$,
some time steps have a negative real part.

Recursively, starting from the second order method $M_{p}$ with $p=2$
and using formula (\ref{eq:real_coefficients}) or (\ref{eq:complex_coefficients})
it is then possible to construct a method of order $p=4$ with $3$
steps, a method of order $p=6$ with $3^{2}=9$ steps, \textit{etc.}

\subsection{Other palindromic schemes}

The above construction allows to constructing symmetric high order
schemes belonging to the family of palindromic schemes \cite{kahan1997composition}.
A general palindromic scheme with $s+1$ steps has the form
\begin{equation}
M_{p}(\Delta t)=M_{2}(\gamma_{0}\Delta t)M_{2}(\gamma_{1}\Delta t)\cdots M_{2}(\gamma_{s}\Delta t),\label{eq:palindromic}
\end{equation}
where the $\gamma_{i}$'s are complex numbers such that
\[
\gamma_{i}=\gamma_{s-i},\quad0\leq i\leq s.
\]
It is possible to find other palindromic schemes with better stability
or precision properties than with formula (\ref{eq:real_coefficients})
or (\ref{eq:complex_coefficients}). For $p=4$ and $s=4$ we have
for example the fourth-order Suzuki scheme (see \cite{suzuki1990fractal,hairer2006geometric,mclachlan2002splitting})
\begin{equation}
\gamma_{0}=\gamma_{1}=\gamma_{3}=\gamma_{4}=\frac{1}{4-4^{1/3}},\quad\gamma_{2}=-\frac{4^{1/3}}{4-4^{1/3}}.\label{eq:suzuki_order4}
\end{equation}

This scheme requires five stages and has one negative time steps but
for all $i$ $\left|\gamma_{i}\right|<1$, which improves stability
and precision compared to the choice (\ref{eq:real_coefficients}).

For $p=6$ and $s=8$, we have the sixth-order Kahan-Li scheme \cite{kahan1997composition}.
This scheme is constructed by writing down the order conditions of
a general nine-stage palindromic method. Among the possible solutions,
the authors of \cite{kahan1997composition} impose the minimization
of some error coefficients. The $\gamma_{i}$'s then have to be computed
numerically. They can be found at \url{http://www.netlib.org/ode/}
and are reproduced below with 60 digits precision:
\begin{equation}
\begin{array}{cc}
\gamma_{0}=\gamma_{8}= & 0.392161444007314139275655330038380932595385404354442882183619\\
\gamma_{1}=\gamma_{7}= & 0.332599136789359438604272125325790569941599549617156528439173\\
\gamma_{2}=\gamma_{6}= & -0.706246172557639359809845337222763994485425050210063375842163\\
\gamma_{3}=\gamma_{5}= & 0.0822135962935508002304427053341134143428469807222103772811280\\
\gamma_{4}= & 0.798543990934829963398950353048958155211186231032507175876486
\end{array}\label{eq:kahan_li_order6}
\end{equation}

In this paper we test some of the previous methods for the orders
$p=2$, $p=4$ and $p=6$.

\section{DG transport solver}

\subsection{Practical implementation}

For the numerical applications,  we need a numerical discretization
$L_{h}$ of the transport operator $L$. As stated above, we have
chosen an upwind high order nodal Discontinuous Galerkin (DG) approximation.
We do not give all the details of the method. We refer for instance
to \cite{hesthaven2007nodal}. We just give a brief description.

First we split the initial interval $[a,b]$ into $N_{x}$ cells of
length $h=(b-a)/N_{x}$. In each cell and for each velocity $v_{i}$
we consider a $d^{th}$ order polynomial approximation of $f_{i}$
in the $x$ variable. The approximation is thus discontinuous at the
cell boundaries. We then apply the upwind DG approximation in order
to construct an approximate transport matrix $L_{h}$ of size $N_{x}\times(d+1)\times N_{v}$
where $N_{v}=2m$. For simplicity reasons, we use a Lagrange polynomial
basis associated to Gauss-Lobatto quadrature points in each cell $[a+jh,a+(j+1)h]$,
$j=0\ldots N_{x}-1$. 

A huge advantage of the upwind scheme is that with a good numbering
of the unknowns, the matrix $L_{h}$ is block triangular, with diagonal
blocks of size $d+1$. Expressed differently, for solving the implicit
upwind DG scheme it is not necessary to assemble and invert a large
linear system, we can compute $f_{i}$ cell after cell, following
the direction of the flow, given by the sign of $v_{i}$. For software
design reasons it is also possible to keep assembly and LU decomposition
procedures. Some linear algebra libraries, such as KLU \cite{davis2010algorithm}
are able to detect in an efficient way block triangular matrices.

In our tests, we apply the composition method $M_{p}(\Delta t)$ with
a time step $\Delta t$ satisfying the relation
\begin{equation}
\Delta t=\beta\frac{\delta}{\lambda},\label{eq:cfl_definition}
\end{equation}
where $\beta$ is the CFL number and $\delta$ is the minimal distance
between two Gauss-Lobatto interpolation points in the cells (recall
that $\lambda$ is the lattice velocity). With an explicit scheme
based on (\ref{eq:good_m2}), we would expect instability for $\beta>\beta_{c}$,
with 
\[
\beta_{c}\simeq\frac{2}{\max_{i}\left|\gamma_{i}\right|}.
\]
In our experiments (see Section \ref{sec:Numerical-results}), we
have verified that in many cases the method remains stable and precise
even for very large values of $\beta.$ However, we have sometimes
observed instabilities, maybe arising from the boundary conditions.
More theoretical investigations are needed in order to understand
the origin of those instabilities.

\subsection{Handling negative time-steps}

As explained in Section \ref{sec:Higher-order-extension} we may have
to apply the elementary collision or transport bricks $C_{2}$ and
$T_{2}$ with a negative time step $-\Delta t<0$. The exact transport
operator is perfectly reversible. If we were using an exact characteristic
solver, negative time steps would not cause any problem. However,
the DG approximation introduces a slight dissipation due to upwinding.
For ensuring stability, we have thus to replace $T_{2}(-\Delta t)$
with a more stable operator. This can be done by first observing that
solving 
\[
\partial_{t}f+v\partial_{x}f=0
\]
for negative time $t<0$ is equivalent to solve
\[
\partial_{t'}f-v\partial_{x}f=0
\]
for $t'=-t>0$. This is a transport equation with the opposite velocity.
Then we also observe that the chosen velocity lattice (\ref{eq:velocity_lattice})
is symmetric. For solving numerically the transport equation with
opposite velocities it is thus enough to exchange at the beginning
and the end of the time step the components of $f$ associated to
$\lambda$ or $-\lambda$. More precisely, we define the exchange
operator
\[
f'=Rf
\]
by
\[
f'_{2k}=f_{2k+1},\quad f'_{2k+1}=f_{2k}.
\]
And in the palindromic scheme, if we encounter a negative time step,
we replace the transport operator
\[
T_{2}(-\Delta t)
\]
by
\[
T'_{2}(-\Delta t)=RT_{2}(\Delta t)R.
\]

\textbf{Remark}: obviously, with this correction, the resulting scheme
is still symmetric in the sense of (\ref{eq:sym_prop}). The time
order property is also preserved if the DG solver is of sufficiently
high order. In practice we have observed that the correction improves
the stability without altering the order of the whole method.

When $\tau=0$ the collision operator $C_{2}$ is perfectly reversible,
as can be seen from (\ref{eq:tau0_sym}), and negative time steps
are not a problem. When $\tau>0$ this can be a problem, especially
if $\Delta t\simeq-2\tau$ (see formula (\ref{eq:crank_nicolson})).
Let us observe, however, that in many cases where $\tau\ll\Delta t$
we have observed stable and precise simulations with the Suzuki (\ref{eq:suzuki_order4})
or Kahan-Li (\ref{eq:kahan_li_order6}) schemes and high CFL numbers.

Theoretical and numerical investigations are still needed to reach
a better understanding of the stability conditions.

\section{Numerical results\label{sec:Numerical-results}}

\subsection{Smooth solution}

In this subsection and the next one we consider an isothermal compressible
flow of density $\rho$ and velocity $u$. The sound speed is fixed
to $c=0.6.$ The conservative system is given by
\[
w=(\rho,\rho u),\quad q(w)=(\rho u,\rho u^{2}+c^{2}\rho).
\]
The eigenvalues of $D_{w}q(w)$ are
\[
\mu_{1}=u-c,\quad\mu_{2}=u+c.
\]

For the first validation of the method we consider a test case with
a smooth solution, in the fluid limit $\tau=0$. This test is performed
with time-independent boundary conditions. Indeed, it is known that
relaxation methods may lead to other difficulties and order reduction
in presence of time-dependent boundary conditions \cite{hundsdorfer2013numerical}.

The initial condition is given by
\[
\rho(x,0)=1+e^{-30x^{2}},\quad u(x,0)=0.
\]
The sound speed is set to $c=0.6$ and the lattice velocity to $\lambda=2$.
For the first convergence study, the CFL number is fixed to $\beta=5$:
see definition (\ref{eq:cfl_definition}). We consider a sufficiently
large computational domain $[a,b]=[-2,2]$ and a sufficiently short
final time $t_{\max}=0.4$ so that the boundary conditions play no
role. The reference solution $f(\cdot,t_{\max})$ is computed numerically
with a very fine mesh and a very small time step. For the composition
methods of order $2$, $4$ and $6$ we adopt a DG solver of order
$6$ in the space variable $x$ (the polynomial order in $x$ is fixed
to $d=5$) for only measuring the time order of the scheme.

\begin{figure}
\centering{}\includegraphics[width=8cm]{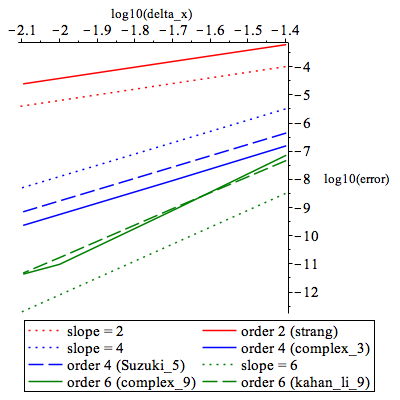}\caption{\label{fig:Convergence-study}Convergence study for several palindromic
methods, order 2 (red), 4 (blue) and 6 (green). The dotted lines are
reference lines with slopes 2, 4 and 6 respectively. At order 4, we
observe that the scheme with complex time steps is more precise than
the five-step Suzuki scheme. At order 6, the scheme with complex time
steps has a similar precision as the Kahan-Li scheme.}
\end{figure}

On Figure \ref{fig:Convergence-study} we give the results of the
convergence study for the smooth solution. The considered error is
the $L^{2}$ norm of $f_{h}(t_{\max})-f(\cdot,t_{\max}).$ 

We observe that the schemes indeed produce the expected accuracy.

\begin{figure}
\centering{}\includegraphics[width=8cm]{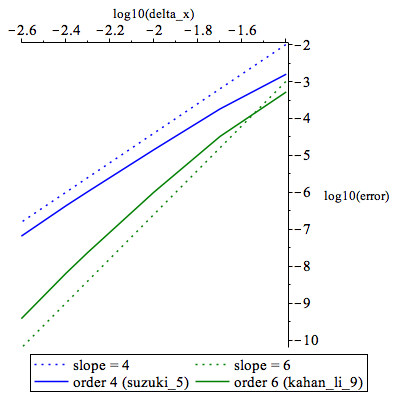}\caption{\label{fig:Convergence-study-beta50}Convergence study for palindromic
methods at order 4 (blue, Suzuki scheme) and 6 (green, Kahan-Li scheme)
with high CFL number $\beta=50$. The dotted lines are reference lines
with slopes 4 and 6 respectively. We observe that the schemes are
stable and that for sufficiently fine meshes, the expected rates of
convergence are attained.}
\end{figure}

We make the same experiment with $\beta=50$. For this choice of CFL
number, we found that the complex time step methods are unstable.
The other methods are stable. The convergence study for the Suzuki
and Kahan-Li schemes is presented on Figure \ref{fig:Convergence-study-beta50}.

\subsection{Behavior for discontinuous solutions}

We have also experimented the scheme with complex time steps for discontinuous
solutions. Of course, in this case the effective order of the method
cannot be higher than one and we expect Gibbs oscillations near the
discontinuities. On the interval $[a,b]=[-1,1],$ we consider a Riemann
problem with the following initial condition
\[
\rho(x,0)=\begin{cases}
2 & \text{if }x<0,\\
1 & \text{otherwise.}
\end{cases},\quad u(x,0)=0.
\]

\begin{figure}
\centering{}\includegraphics[width=15cm]{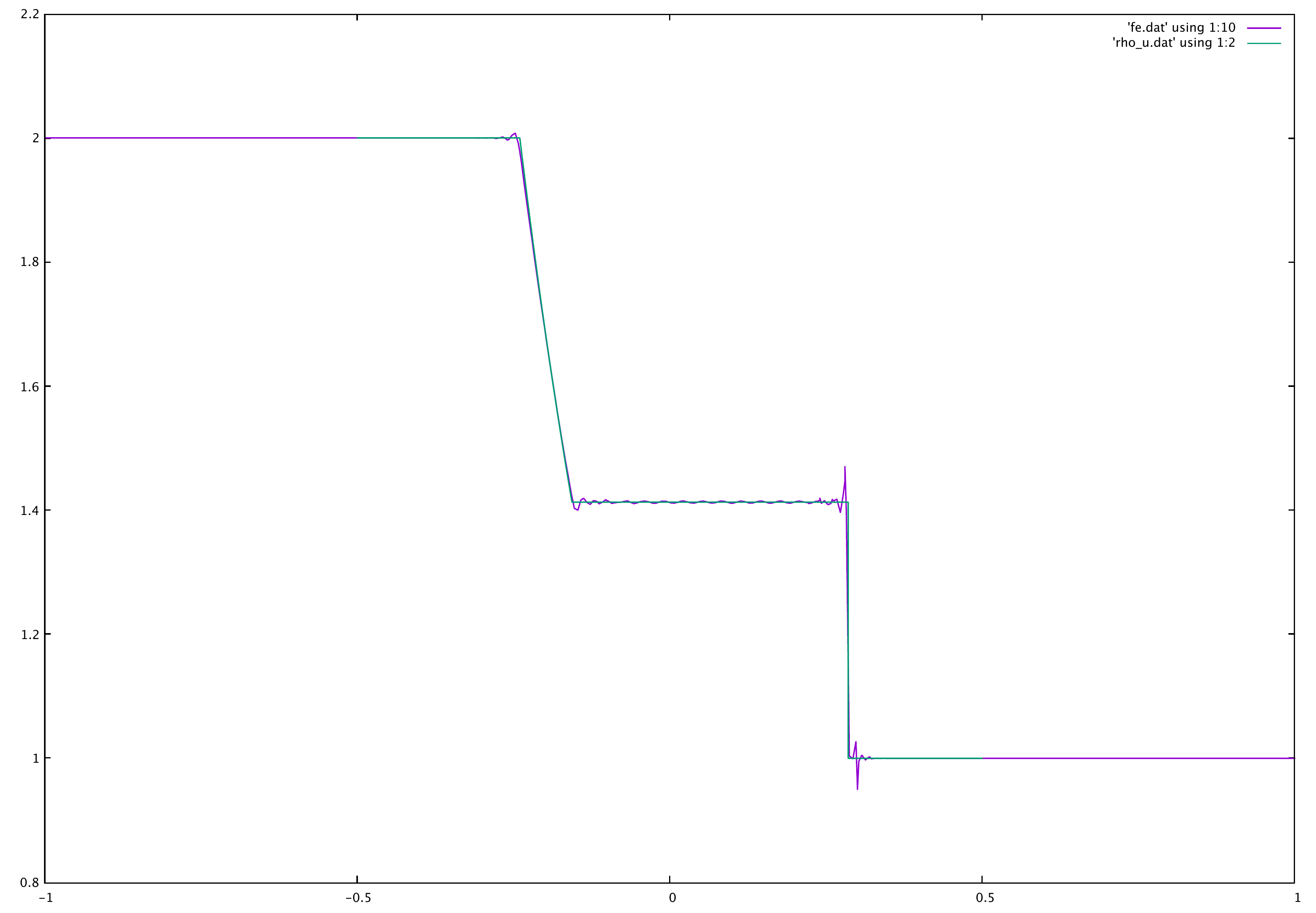}\caption{\label{fig:riemann_rho}Riemann problem with $\tau=0$. Comparison
of the exact density (green curve), and the numerical sixth-order
solution with complex time steps (purple curve). The number of DG
cells is $N_{x}=100$ and the CFL number $\beta=3$.}
\end{figure}

We consider numerical results in the fluid limit $\tau=0$. On Figure
\ref{fig:riemann_rho} we compare the sixth-order numerical solution
with the exact one for a CFL number $\beta=3$ and $N_{x}=100$ cells.
We observe that the high order scheme is able to capture a precise
rarefaction wave and the correct position of the shock wave. We observe
oscillations in the shock wave as expected.

\begin{figure}
\centering{}\includegraphics[width=15cm]{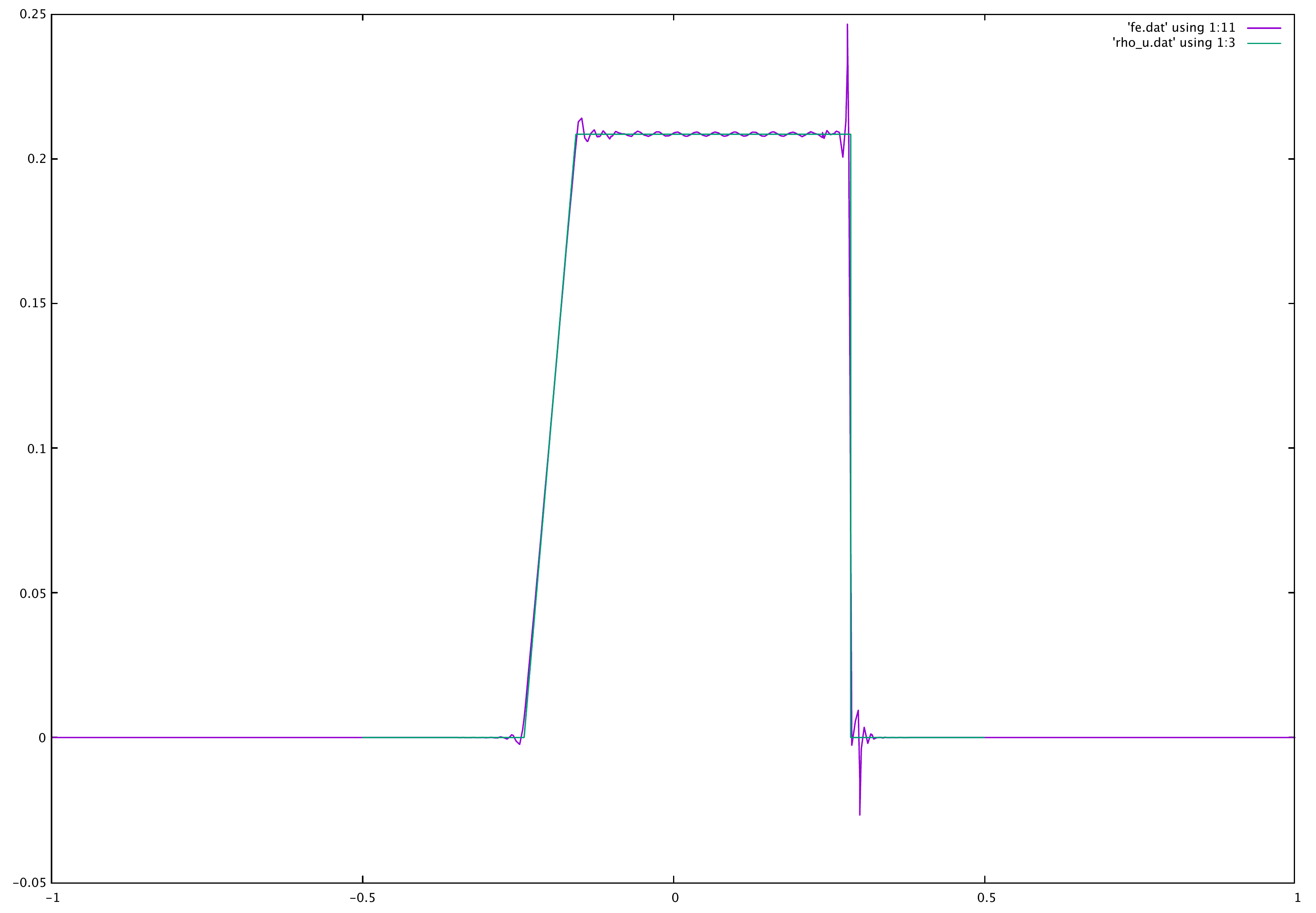}\caption{\label{fig:riemann_u}Riemann problem with $\tau=0$. Comparison of
the exact velocity (green curve), and the numerical sixth-order solution
with complex time steps (purple curve). The number of DG cells is
$N_{x}=100$ and the CFL number $\beta=3$.}
\end{figure}
\begin{figure}
\centering{}\includegraphics[width=15cm]{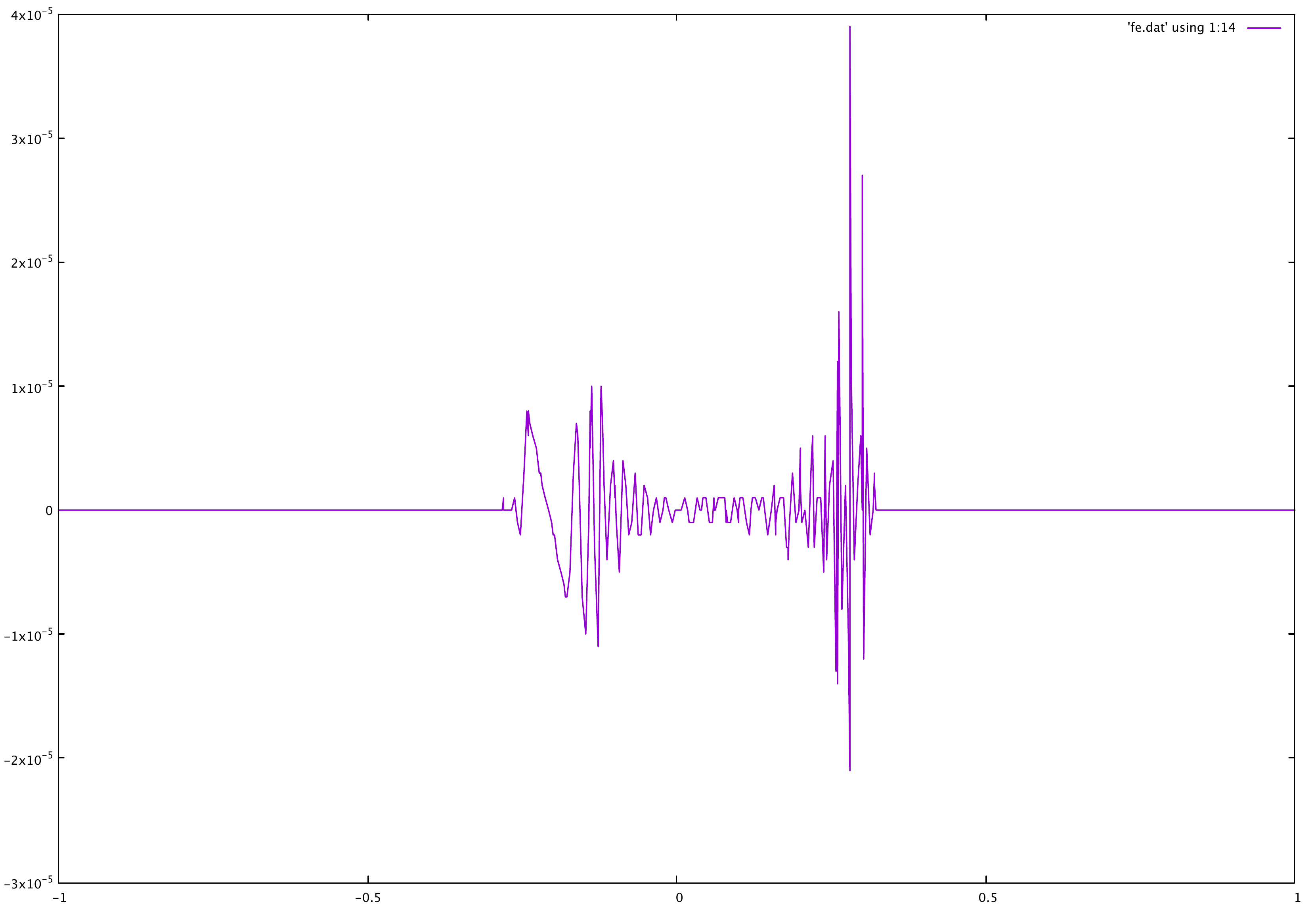}\caption{\label{fig:imaginary_rho}Riemann problem with $\tau=0$. Imaginary
part of the numerical density. We observe that it is of order $10^{-5}$
and that it has higher variations near to the shock. Maybe that it
could be used as an error indicator.}
\end{figure}

The complex time step method is used for approximating the real solution
of a differential equation. We thus expect that the imaginary part
of the numerical solution is small. We plot the imaginary part of
$\rho$ on Figure \ref{fig:imaginary_rho}. We observe that the imaginary
part is indeed small (of the order $10^{-5}$) and that it takes its
higher values in the discontinuous region. Maybe that it could be
used as a shock indicator for detecting oscillations and controlling
them. 

\subsection{Low Mach flow}

\begin{figure}
\centering{}\includegraphics[width=8cm]{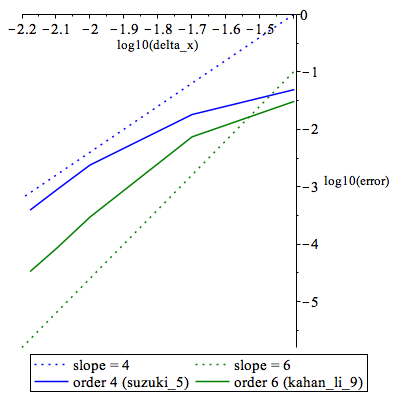}\caption{\label{fig:Convergence-study-lowmach}Convergence study for the low
Mach test case with a CFL number $\beta=100$. Comparison of the fourth-order
Suzuki scheme (blue) and the sixth-order Kahan-Li scheme (green).
The dotted lines are reference lines with slopes 4 and 6 respectively.
We observe the expected convergence rates for fine enough meshes.
With this CFL number, the schemes with complex time steps are unstable.}
\end{figure}

In this section, we consider the Euler equations for a compressible
gas with polytropic exponent $\gamma=7/5$. The primitive unknowns
are the density $\rho$, the velocity $u$ and the pressure $p$.
The conservative system is given by
\begin{equation}
w=(\rho,\rho u,\frac{1}{2}\rho u^{2}+\frac{p}{\gamma-1}),\quad q(w)=(\rho u,\rho u^{2}+p,\frac{1}{2}\rho u^{3}+\frac{\gamma pu}{\gamma-1}).\label{eq:euler_energy}
\end{equation}
The eigenvalues of $D_{w}q(w)$ are
\[
\mu_{1}=u-c,\quad\mu_{2}=u,\quad\mu_{3}=u+c,\quad\text{ with }c=\sqrt{\frac{\gamma p}{\rho}}.
\]
We consider the following exact solution
\[
u=\frac{1}{100},\quad p=1,\quad\rho=\omega\rho_{L}+(1-\omega)\rho_{R},
\]
with
\[
\rho_{L}=2,\quad\rho_{R}=1,\quad\omega=\left(1-\text{erf}(10(x-ut))\right)/2.
\]
It represents a smooth and slowly moving contact wave. We compare
the exact and numerical solutions at time $t=20$ with a CFL number
$\beta=100$ and $\lambda=2$. The convergence study is presented
on Figure \ref{fig:Convergence-study-lowmach}. We observe the expected
convergence rates. For such CFL number $\beta$, the schemes with
complex time steps are unstable.

\subsection{Viscous test}

In the above tests, we have assumed infinitely fast relaxation $\tau=0$.
In this subsection we present a test with $\tau>0$. In practice the
relaxation schemes are generally used with small relaxation time $\tau$.
Indeed, with Taylor or Chapman-Enskog expansions, it is possible to
prove that the relaxation scheme is, at order $O(\tau^{2})$, an approximation
of
\[
\partial_{t}w+\partial_{x}q(w)-\tau\partial_{x}B\partial_{x}w=0,
\]
where the matrix $B$ depends on the choice of the lattice velocities,
Maxwellian states, \textit{etc}. To some extent it is possible to
adjust $B$ to real physical terms, such as the Navier-Stokes viscosity
(on this subject, see for instance \cite{dubois2016recovering,qian1992lattice,he1997lattice}).
The higher order terms in $\tau$ involve higher order space derivatives,
which are generally undesirable. The smallness of $\tau$ implies
that, in practice, it is important to handle the cases $\tau<\Delta t$
or $\tau\approx\Delta t$. 

We consider the Euler model (\ref{eq:euler_energy}) with $x\in[-1/2,1/2]$,
$t\in[0,0.2]$ and the following initial condition
\[
w(x,0)=\begin{cases}
w_{L} & \text{if }x<0,\\
w_{R} & \text{otherwise.}
\end{cases}
\]
The left and right states are given by
\[
w_{L}=(2,0,5),\quad w_{R}=(1,0,5/2).
\]
We consider the sixth-order Kahan-Li scheme (\ref{eq:kahan_li_order6})
with a mesh of size $N_{x}=100$ and a CFL number $\beta=5$.

We have already observed that the collision step (\ref{eq:crank_nicolson})
cannot be applied if $2\tau+\gamma_{2}\Delta t/2=0$ (the $\Delta t/2$
comes from the definition of the elementary brick (\ref{eq:good_m2})).
In this way we find a critical value of $\tau=\tau_{c}=0.000519$
for which the Kahan-Li scheme is unusable. With the same parameters,
we have found that the five-step Suzuki scheme (\ref{eq:suzuki_order4})
is unstable.

\begin{figure}
\centering{}\includegraphics[width=15cm]{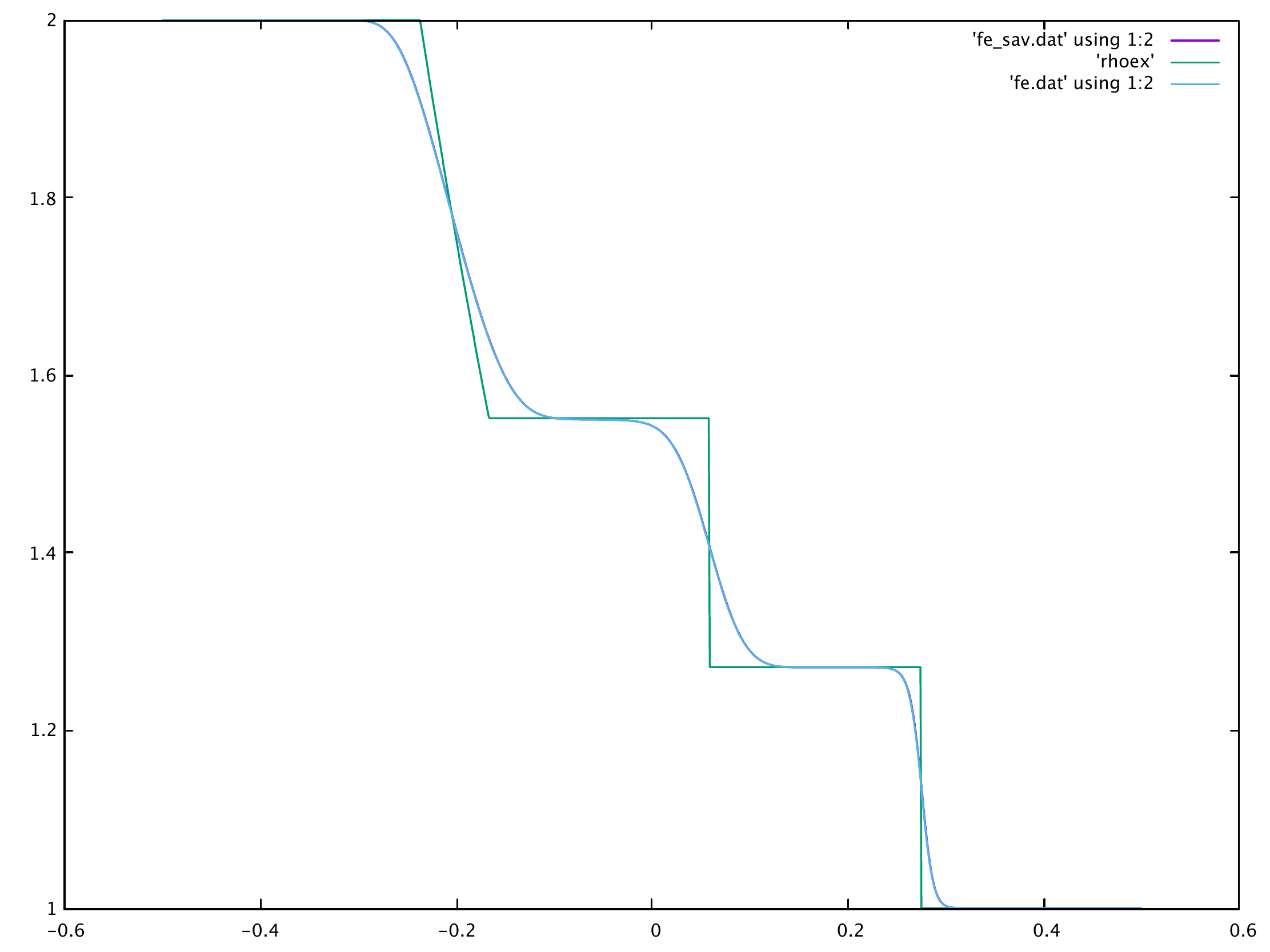}\caption{\label{fig:riemann_rho_viscous}Riemann problem with $\tau=0.000519$.
Comparison of the exact density of the inviscid case (green curve),
and the sixth-order complex numerical solution with $\tau>0$ (blue
and purple curves: indistinguishable). The number of DG cells is $N_{x}=100$
or $N_{x}=1000$ and the CFL number $\beta=5$. The relaxation parameter
has the same effect as a viscosity parameter.}
\end{figure}

However, the sixth-order scheme with complex time steps (\ref{eq:complex_coefficients})
is stable. The resulting density is plot on Figure \ref{fig:riemann_rho_viscous}.
On the same figure we also represent a computation on a finer mesh
with $N_{x}=1000$ cells. The curves with $N_{x}=100$ and $N_{x}=1000$
are indistinguishable (the error in the $L^{2}$ norm is of the order
of $10^{-4}).$

If we take a higher CFL number $\beta=10,$ we then observe that the
Kahan-Li scheme becomes stable again and produces results that are
very similar to those of Figure \ref{fig:riemann_rho_viscous}.

The conclusion of this numerical test is that it is possible to compute
with high precision the solution of the relaxation system when $\tau>0$.
However, the stability conditions of the different schemes are for
the moment unclear. It seems that it is possible to cross the ``$2\tau+\Delta t$
barrier'' if we consider the complex method at moderate CFL numbers
$\beta.$ 

\section{Conclusion}

In this paper we have described a new numerical method, the Palindromic
Discontinuous Galerkin Lattice Boltzmann Method, for solving kinetic
equations with stiff relaxation. The method has the following features:
\begin{itemize}
\item Time integration is high order, based on a general palindromic composition
method. We have tested it at order $4$ and $6$.
\item The transport solver is based on an implicit high order upwind Discontinuous
Galerkin method. Thanks to the upwind flux, the linear system to be
solved at each time step is triangular.
\item The scheme has better stability properties than an explicit scheme
and allows for low storage optimization.
\item The method is general, highly parallel and can be extended to higher
dimensions.
\end{itemize}
We are currently working on the extension of the method to higher
dimensions and to optimizations of the implementation on hybrid computers.
A lot of theoretical and numerical investigations are still needed
for understanding the stability properties of the method. It is also
important in practical applications to extend the method with more
general boundary conditions. 

\section*{Bibliography}

\bibliographystyle{plain}
\bibliography{helluy-relax}

\begin{thebibliography}{10}

\bibitem{bey1997downwind}
J{\"u}rgen Bey and Gabriel Wittum.
\newblock Downwind numbering: Robust multigrid for convection-diffusion
  problems.
\newblock {\em Applied Numerical Mathematics}, 23(1):177--192, 1997.

\bibitem{castella2009splitting}
Fran{\c{c}}ois Castella, Philippe Chartier, St{\'e}phane Descombes, and Gilles
  Vilmart.
\newblock Splitting methods with complex times for parabolic equations.
\newblock {\em BIT Numerical Mathematics}, 49(3):487--508, 2009.

\bibitem{chen1998lattice}
Shiyi Chen and Gary~D Doolen.
\newblock Lattice {B}oltzmann method for fluid flows.
\newblock {\em Annual review of fluid mechanics}, 30(1):329--364, 1998.

\bibitem{coquel2008large}
Fr{\'e}d{\'e}ric Coquel, Q-L Nguyen, Marie Postel, and Q-H Tran.
\newblock Large time step positivity-preserving method for multiphase flows.
\newblock In {\em Hyperbolic Problems: Theory, Numerics, Applications}, pages
  849--856. Springer, 2008.

\bibitem{davis2010algorithm}
Timothy~A Davis and Ekanathan Palamadai~Natarajan.
\newblock Algorithm 907: {KLU}, a direct sparse solver for circuit simulation
  problems.
\newblock {\em ACM Transactions on Mathematical Software (TOMS)}, 37(3):36,
  2010.

\bibitem{dellar2013interpretation}
Paul~J Dellar.
\newblock An interpretation and derivation of the lattice {B}oltzmann method
  using {S}trang splitting.
\newblock {\em Computers \& Mathematics with Applications}, 65(2):129--141,
  2013.

\bibitem{dubois2013stable}
Fran{\c{c}}ois Dubois.
\newblock Stable lattice {B}oltzmann schemes with a dual entropy approach for
  monodimensional nonlinear waves.
\newblock {\em Computers \& Mathematics with Applications}, 65(2):142--159,
  2013.

\bibitem{dubois2016recovering}
Fran{\c{c}}ois Dubois, Benjamin Graille, and Pierre Lallemand.
\newblock Recovering the full {N}avier-{S}tokes equations with lattice
  {B}oltzmann schemes.
\newblock In {\em 30th international symposium on rarefied gas dynamics: rgd
  30}, volume 1786, page 040003. AIP Publishing, 2016.

\bibitem{graille2014approximation}
Benjamin Graille.
\newblock Approximation of mono-dimensional hyperbolic systems: A lattice
  {B}oltzmann scheme as a relaxation method.
\newblock {\em Journal of Computational Physics}, 266:74--88, 2014.

\bibitem{hairer2006geometric}
Ernst Hairer, Christian Lubich, and Gerhard Wanner.
\newblock {\em Geometric numerical integration: structure-preserving algorithms
  for ordinary differential equations}, volume~31.
\newblock Springer Science \& Business Media, 2006.

\bibitem{he1997lattice}
Xiaoyi He and Li-Shi Luo.
\newblock Lattice {B}oltzmann model for the incompressible {N}avier--{S}tokes
  equation.
\newblock {\em Journal of statistical Physics}, 88(3-4):927--944, 1997.

\bibitem{helluy:hal-01403759}
Philippe Helluy.
\newblock {Stability analysis of an implicit lattice {B}oltzmann scheme}.
\newblock To appear in Oberwolfach reports.
  https://hal.archives-ouvertes.fr/hal-01403759, 2016.

\bibitem{hesthaven2007nodal}
Jan~S Hesthaven and Tim Warburton.
\newblock {\em Nodal discontinuous Galerkin methods: algorithms, analysis, and
  applications}.
\newblock Springer Science \& Business Media, 2007.

\bibitem{hundsdorfer2013numerical}
Willem Hundsdorfer and Jan~G Verwer.
\newblock {\em Numerical solution of time-dependent
  advection-diffusion-reaction equations}, volume~33.
\newblock Springer Science \& Business Media, 2013.

\bibitem{jin1995runge}
Shi Jin.
\newblock Runge-{K}utta methods for hyperbolic conservation laws with stiff
  relaxation terms.
\newblock {\em Journal of Computational Physics}, 122(1):51--67, 1995.

\bibitem{jin1999efficient}
Shi Jin.
\newblock Efficient asymptotic-preserving {(AP)} schemes for some multiscale
  kinetic equations.
\newblock {\em SIAM Journal on Scientific Computing}, 21(2):441--454, 1999.

\bibitem{jin1995relaxation}
Shi Jin and Zhouping Xin.
\newblock The relaxation schemes for systems of conservation laws in arbitrary
  space dimensions.
\newblock {\em Communications on pure and applied mathematics}, 48(3):235--276,
  1995.

\bibitem{kahan1997composition}
William Kahan and Ren-Cang Li.
\newblock Composition constants for raising the orders of unconventional
  schemes for ordinary differential equations.
\newblock {\em Mathematics of Computation of the American Mathematical
  Society}, 66(219):1089--1099, 1997.

\bibitem{liu1987hyperbolic}
Tai-Ping Liu.
\newblock Hyperbolic conservation laws with relaxation.
\newblock {\em Communications in Mathematical Physics}, 108(1):153--175, 1987.

\bibitem{mclachlan2002splitting}
Robert~I McLachlan and G~Reinout~W Quispel.
\newblock Splitting methods.
\newblock {\em Acta Numerica}, 11:341--434, 2002.

\bibitem{mei1998finite}
Renwei Mei and Wei Shyy.
\newblock On the finite difference-based lattice {B}oltzmann method in
  curvilinear coordinates.
\newblock {\em Journal of Computational Physics}, 143(2):426--448, 1998.

\bibitem{min2011spectral}
Misun Min and Taehun Lee.
\newblock A spectral-element {D}iscontinuous {G}alerkin lattice {B}oltzmann
  method for nearly incompressible flows.
\newblock {\em Journal of Computational Physics}, 230(1):245--259, 2011.

\bibitem{moustafa20153d}
Salli Moustafa, Mathieu Faverge, Laurent Plagne, and Pierre Ramet.
\newblock {3D} cartesian transport sweep for massively parallel architectures
  with {PARSEC}.
\newblock In {\em Parallel and Distributed Processing Symposium (IPDPS), 2015
  IEEE International}, pages 581--590. IEEE, 2015.

\bibitem{nannelli1992lattice}
Francesca Nannelli and Sauro Succi.
\newblock The lattice {B}oltzmann equation on irregular lattices.
\newblock {\em Journal of Statistical Physics}, 68(3-4):401--407, 1992.

\bibitem{natvig2008fast}
Jostein~R Natvig and Knut-Andreas Lie.
\newblock Fast computation of multiphase flow in porous media by implicit
  discontinuous galerkin schemes with optimal ordering of elements.
\newblock {\em Journal of Computational Physics}, 227(24):10108--10124, 2008.

\bibitem{peng1998lattice}
Gongwen Peng, Haowen Xi, Comer Duncan, and So-Hsiang Chou.
\newblock Lattice {B}oltzmann method on irregular meshes.
\newblock {\em Physical Review E}, 58(4):R4124, 1998.

\bibitem{qian1992lattice}
YH~Qian, Dominique d'Humi{\`e}res, and Pierre Lallemand.
\newblock Lattice {BGK} models for {N}avier-{S}tokes equation.
\newblock {\em EPL (Europhysics Letters)}, 17(6):479, 1992.

\bibitem{shi2003discontinuous}
Xing Shi, Jianzhong Lin, and Zhaosheng Yu.
\newblock {D}iscontinuous {G}alerkin spectral element lattice {B}oltzmann
  method on triangular element.
\newblock {\em International Journal for Numerical Methods in Fluids},
  42(11):1249--1261, 2003.

\bibitem{suzuki1990fractal}
Masuo Suzuki.
\newblock Fractal decomposition of exponential operators with applications to
  many-body theories and monte carlo simulations.
\newblock {\em Physics Letters A}, 146(6):319--323, 1990.

\bibitem{wang1999crosswind}
Feng Wang and Jinchao Xu.
\newblock A crosswind block iterative method for convection-dominated problems.
\newblock {\em SIAM Journal on Scientific Computing}, 21(2):620--645, 1999.

\bibitem{witham1974linear}
Gerald~B Witham.
\newblock {\em Linear and Non-linear Waves}.
\newblock Wiley, 1974.

\end{thebibliography}

\end{document}